\documentclass[12pt]{elsarticle}

\usepackage{amssymb}
\usepackage{amsthm}

\usepackage{amsmath,amsfonts,color,verbatim}
\usepackage{amscd}
\usepackage[config,labelfont={rm,rm},textfont=rm]{caption,subfig}
\usepackage{bm} 
\usepackage{wrapfig}
\newcommand{\curl}{\ensuremath{\operatorname{\mathbf{curl}}}}
\newcommand{\grad}{\ensuremath{\operatorname{\mathbf{grad}}}}
\newcommand{\rot}{\ensuremath{\operatorname{rot}}}
\newcommand{\ddiv}{\ensuremath{\operatorname{div}}}
\newcommand{\Hcurl}{\ensuremath{\mathbf{H}(\curl)}}
\newcommand{\Hdiv}{\ensuremath{H(\ddiv)}}

\newcommand{\meas}{\ensuremath{\operatorname{meas}}}
\newcommand{\diag}{\ensuremath{\operatorname{diag}}}
\newtheorem{theorem}{Theorem}

\theoremstyle{remark}
\newtheorem{remark}[theorem]{Remark}

\numberwithin{equation}{section}
\journal{Computers \& Mathematics with Applications}

\graphicspath{{./}{./}} 

\numberwithin{figure}{section}

\begin{document}

\begin{frontmatter}



\title{A Finite Element Framework for Some Mimetic Finite Difference Discretizations}


\author[uz]{C.~Rodrigo\corref{cor1}}
\author[uz]{F.J.~Gaspar}
\author[tufts]{X.~Hu}
\author[penn,sof]{L.~Zikatanov}
\cortext[cor1]{Corresponding author. Phone: +34 976761142; E-mail: carmenr@unizar.es}
\address[uz]{Applied Mathematics Department and IUMA, University of Zaragoza, Spain.}
\address[tufts]{Department of Mathematics, Tufts University, Medford, MA, USA.}
\address[penn]{Department of Mathematics, Penn State University, University Park,
PA, USA}
\address[sof]{Institute for Mathematics and Informatics, Bulgarian
Academy of Sciences, Sofia, Bulgaria}

\begin{abstract}
  In this work we derive equivalence relations between mimetic finite
  difference schemes on simplicial grids and modified
  N\'ed\'elec-Raviart-Thomas finite element methods for model problems
  in $\mathbf{H}(\operatorname{\mathbf{curl}})$ and $H(\operatorname{div})$. This provides a simple and transparent way to
  analyze such mimetic finite difference discretizations using
  the well-known results from finite element theory. The finite element
  framework that we develop is also crucial for the design of efficient
  multigrid methods for mimetic finite difference discretizations,
  since it allows us to use canonical inter-grid transfer operators
  arising from the finite element framework. We provide special Local
  Fourier Analysis and numerical results to demonstrate the
  efficiency of such multigrid methods.
\end{abstract}


\end{frontmatter}




\section{Introduction}\label{sec:intro}


We consider mimetic finite difference (MFD) methods for problems in
$\Hcurl$ and $\Hdiv$ with essential
boundary conditions. Such methods are designed in order to have
natural discrete analogues of conservation (of mass, momentum, etc),
symmetry and positivity of the operators. They are also structure
preserving discretizations, namely, they form discrete de~Rham
complexes.

Such discretization techniques were started in the School of
A.~A.~Samarskii at the Moscow State University, and they have been
further developed and analyzed by Shashkov~\cite{Shashkov_Book} and
Vabishchevich~\cite{Mimetic_Vabishchevich}.
Regarding the MFD methods, our presentation here
follows Vabishchevich~\cite{Mimetic_Vabishchevich} and his
Vector Analysis Grid Operators (VAGO) framework for dual
simplicial/polyhedral (Delaunay/Voronoi) grids.

Many authors have contributed to the research in this field, by applying
the MFD methods successfully to several applications ranging
from diffusion~\cite{Hyman1997, Morel1998, Shashkov1996}, magnetic
diffusion and electromagnetics~\cite{Hyman2001} to continuum
mechanics~\cite{Margolin2000} and gas dynamics~\cite{Campbell2001}. We
refer to a recent comprehensive review paper by Lipnikov, Manzini, and
Shashkov~\cite{Lipnikov2014} and a recent book by Beir\~{a}o da Veiga, Lipnikov, and Manzini~\cite{Veiga.L;Lipnikov.K;Manzini.G2014a} on the subject for details and literature
review.

We are interested in the MFD discretizations of two (standard) model
problems in \Hcurl \; and \Hdiv.  We show that the MFD methods can be
fitted in a more or less standard finite element (FE) framework which leads
to convergence results and makes the design of efficient and fast
solvers for the resulting linear systems quite easy.  Our approach is somewhat like special discrete Hodge operates and, therefore, is related to the generalized finite difference approach proposed by Bossavit (see e.g.~\cite{Bossavit.A2005a} and references therein).  We point out that, in the classical finite difference setting, convergence results
exist, as can be seen in~\cite{Mimetic_Vabishchevich}, but deriving
them is by all means not an easy task.   Moreover, while we
provide details on the constructions in 2D, the equivalence between
the MFD methods and the FE methods carries over
with trivial modifications to 3D case as well. We have only chosen 2D
because it makes the exposition much easier to understand.


Such connections between the MFD schemes and the
mixed FE methods for diffusion equations
with Raviart­-Thomas elements have been already established,
see~\cite{Berndt2001,Berndt2005,Brezzi2005,Kuznetsov2003,Kuznetsov2004}
and references therein. In fact, designing finite element methods on
arbitrary grids is a hot topic and we refer to the recent works on
agglomerated
grids~\cite{Lashuk_Vassilevski2012,Lashuk_Vassilevski2014,Pasciak_Vassilevski} and virtual
finite element methods~\cite{Beirao2013,Beirao2013-hitch,Brezzi2014}.

Most of the existing works are on approximation, stability and
structure preserving properties of the MFD discretizations. Developing fast
solvers for the resulting linear systems is a topic that needs more
attention, since the design of fast solvers makes the MFD
discretizations more practical and efficient. For FE methods, solvers can be built using the
agglomeration techniques introduced by Lashuk and
Vassilevski~\cite{Lashuk_Vassilevski2012,Lashuk_Vassilevski2014}. Such
techniques do not apply to the MFD discretizations
(even on simplicial grids!) and, to the best of our knowledge, such
results are not available in the literature.  We point out though that on rectangular grids for standard
finite difference schemes for $\Hdiv$ problems, a distributive relaxation based multigrid was
proposed in~\cite{Gaspar_grad_div}.

As we have pointed out, our goal is to apply classical multigrid and
subspace correction techniques~\cite{Brandt77, Hackb, TOS01, Wess,
Jinchao92} for the mimetic discretizations, by first establishing the
relation with N\'ed\'elec-Raviart-Thomas elements.  Such
approach automatically makes efficient methods such as the ones
developed by Arnold, Falk and Winther~\cite{Arnold_et_al} and Hiptmair
and Xu (HX)~\cite{Hiptmair2007} preconditioners applicable for the MFD methods.

Regarding the convergence of $W-$ and $V-$cycle
multigrid with a multiplicative Schwarz relaxation proposed
in~\cite{Arnold_et_al}, we complement the numerical results with
practical Local Fourier Analysis
(LFA) which provides sharp estimates of the multigrid convergence rates. We use a variant
of LFA that is applicable on simplicial grids (see~\cite{RodrigoBook}) and
compare the convergence rates predicted by LFA with the
actual convergence rates of $W-$cycle and $V-$cycle multigrid.

The rest of the paper is organized as follows. In
Section~\ref{sec:mimeticFD}, we describe the MFD
schemes on simplicial grids.
In Section~\ref{sec:FD_FEM} we derive the ``modified''
N\'ed\'elec-Raviart-Thomas FE methods and show their
equivalence to the VAGO MFD schemes.
Section~\ref{sec:FEM_mg} defines the
multigrid components: smoothers, and, with the help of the results
from Section~\ref{sec:FD_FEM}, the canonical inter-grid transfer
operators. In this section, we also discuss the setup and the design
of appropriate LFA for edge-based discretizations and Schwarz smoothers.
The results obtained from the LFA analysis are shown in
Section~\ref{sec:experiments}, together with the convergence rates
of the resulting multigrid algorithm. Finally, conclusions are drawn in
Section~\ref{sec:conclusions}.

\section{Mimetic finite difference discretizations on triangular grids}\label{sec:mimeticFD}

We consider  the following two model problems for $\mathbf u$
in a two dimensional
simply connected domain $\Omega$:
\begin{eqnarray}
&& \curl\rot \mathbf  u + \kappa \mathbf  u = \mathbf   f, \quad \mbox{in} \; \Omega, \label{rot_rot_problem}\\
&& -\grad \ddiv \mathbf u + \kappa \mathbf  u = \mathbf  f, \quad \mbox{in} \; \Omega, \label{grad_div_problem}
\end{eqnarray}
subject to essential boundary conditions (vanishing tangential or
normal components respectively). We also used $\mathbf u$ and
$\mathbf f$ to denote solutions and right hand sides for both
problems without distinguish them in \emph{different} equations and spaces explicitly.
The corresponding variational forms (used in the derivation
of the FE scheme) are: Find ${\mathbf u}\in \Hcurl$ and ${\mathbf u}\in \Hdiv$, respectively, such that
\begin{eqnarray}
&&(\rot \mathbf u , \rot \mathbf  v) + \kappa (\mathbf u,\mathbf v) =
   (\mathbf f, \mathbf v), \quad \mbox{for all} \
   {\mathbf v}\in \Hcurl,\label{rot_rot_problem-w}\\
&&(\ddiv \mathbf u , \ddiv \mathbf v) + \kappa (\mathbf u,\mathbf v) =
   (\mathbf f,\mathbf v), \quad \mbox{for all} \
   {\mathbf v}\in \Hdiv. \label{grad_div_problem-w}
\end{eqnarray}
In 3D we replace $\rot$ with a 3-dimensional $\curl$. In the
variational form, $\Hcurl$ and $\Hdiv$, are the
spaces of square integrable vector valued functions which also have
square integrable $\rot$ ($\curl$ in 3D) or $\ddiv$ respectively. The
functions in the spaces $\Hcurl$ and $\Hdiv$ are also assumed to
satisfy the essential boundary conditions
$(\mathbf{u}\times \mathbf  n)=0$ for~\eqref{rot_rot_problem-w} and
$(\mathbf{u}\cdot \mathbf  n)=0$ for~\eqref{grad_div_problem-w} where $\mathbf{n}$ is the unit normal vector outward to $\partial \Omega$.

\subsection{Mimetic finite differences on a pair of dual meshes}
\begin{wrapfigure}{r}{0.5\textwidth}
  \begin{center}
\includegraphics*[width=0.48\textwidth]{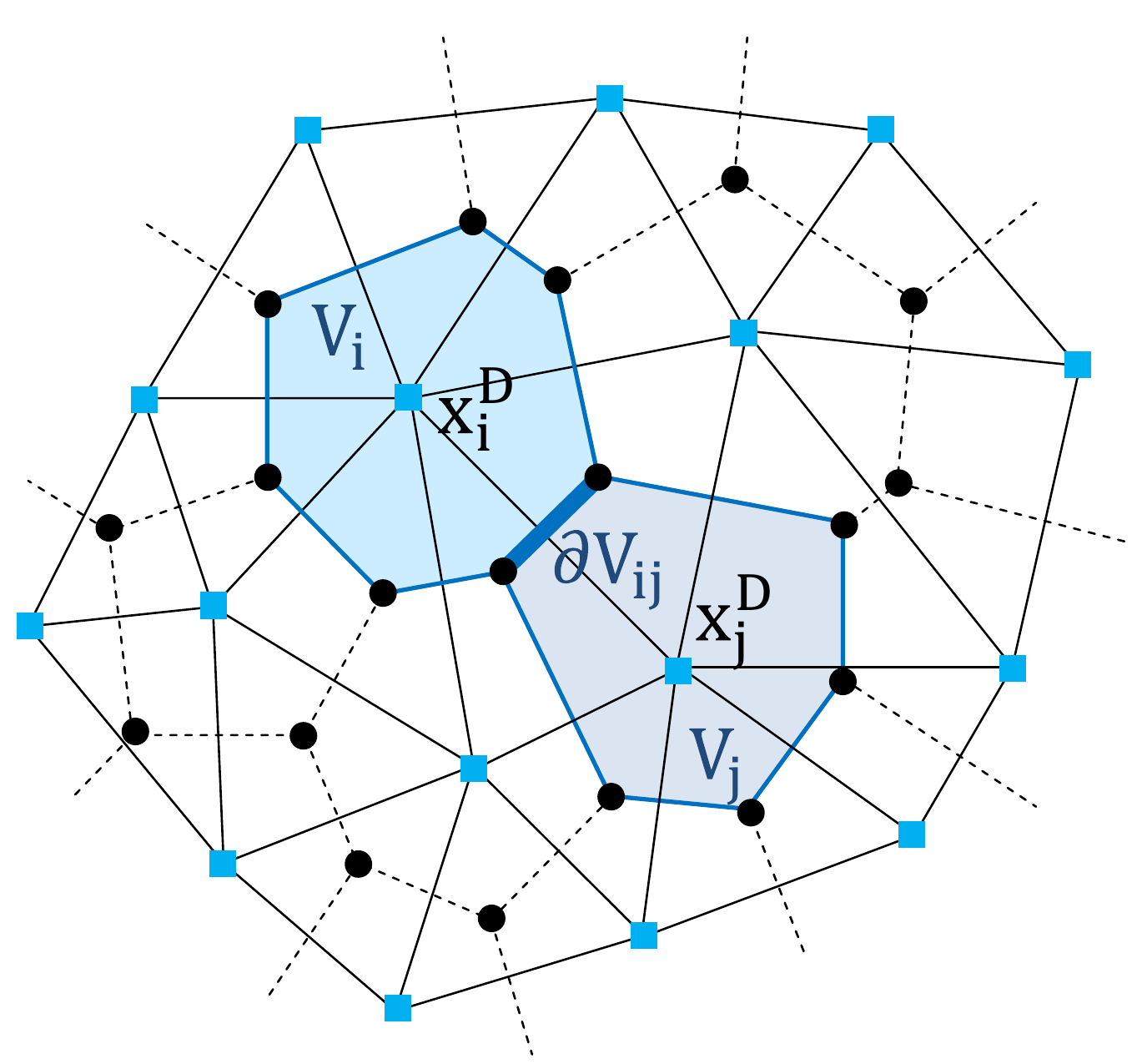}
\end{center}
\caption{A Delaunay mesh and its dual Voronoi mesh.\label{Delaunay_Voronoi_grids}}
\end{wrapfigure}
 We consider MFD schemes for~\eqref{rot_rot_problem} and~\eqref{grad_div_problem} discretized on a pair of primal (Delaunay) simplicial grid and a dual (Voronoi) polyhedral grid.
The vertices of the Delaunay triangulation are
$\{\mathbf{x}_i^D\}_{i=1}^{N_D},$ and the vertices of its dual
Voronoi mesh are the circumcenter of the
Delaunay triangles. We denote the Voronoi vertices by
$\{\mathbf{x}_k^V\}_{k=1}^{N_V}$, and note that each such vertex corresponds to a
Delaunay triangle $D_k$, for $k=1,\ldots,N_V$. In
Figure~\ref{Delaunay_Voronoi_grids} we have depicted a pair of dual
meshes and marked the
Delaunay grid-points by squares and the Voronoi grid-points by circles.  As is
typical in the MFD schemes, we assume that all triangles in the triangulation have
only acute angles. This assumption guarantees that the Voronoi
vertices will always be in the interior of the Delaunay triangles. For
3D analogues of this assumption we refer
to~\cite{Mimetic_Vabishchevich}.
 By duality, to  a Delaunay grid point ${\mathbf x}_i^D$, there
corresponds a Voronoi polygon $V_i$,
\[
V_i = \{{\mathbf x}\in\Omega \, \vert \, \vert {\mathbf x}-{\mathbf
  x}_i^D \vert \le \vert {\mathbf x}-{\mathbf x}_j^D \vert, \,
j=1,\ldots,N_D, \, j\neq i \},
\]
and we denote the Voronoi edge $ V_{ij} = \partial V_i \cap \partial V_j$.

\begin{figure}[htb] \centering
\subfloat[Orientation of the Delaunay edges.]
{
\includegraphics*[width = 4cm]{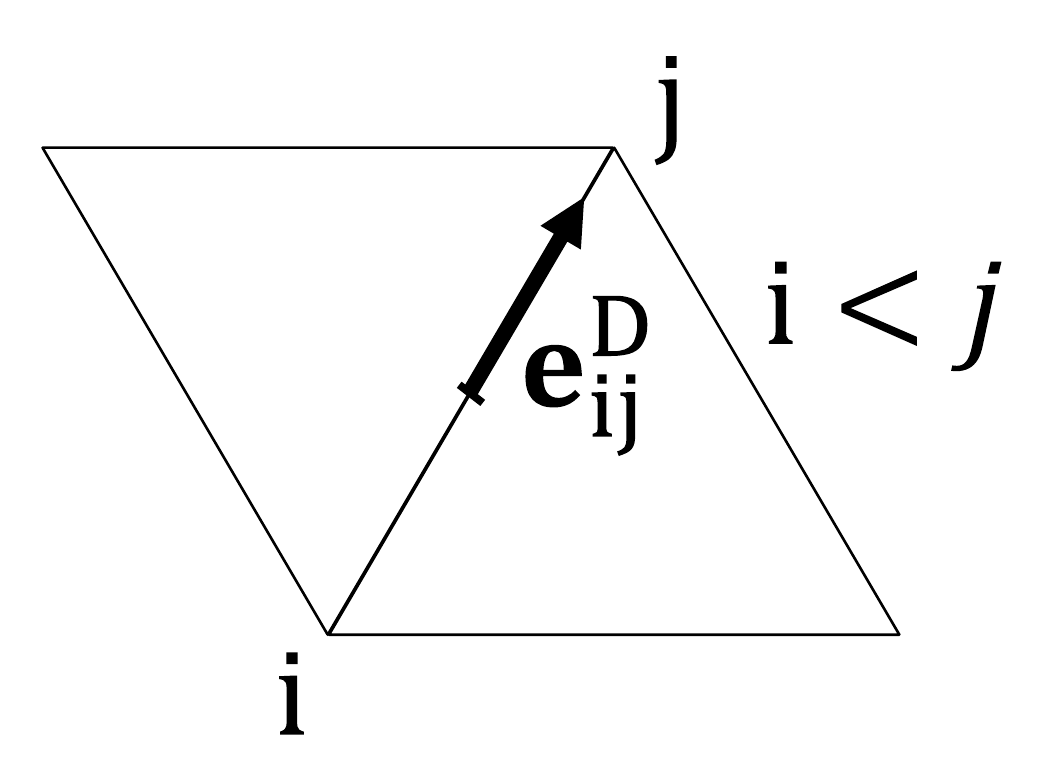}
\label{de-a}
}\hspace*{2em}
\subfloat[Voronoi polygon and notation for the divergence operator.]
{
\includegraphics*[scale = 0.5]{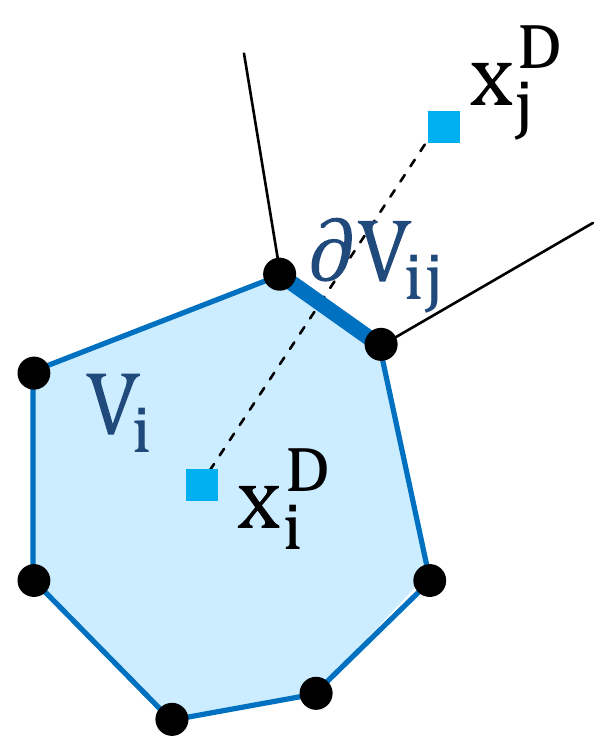}
\label{de-b}
}
\caption{}\label{direction_edges}
\end{figure}

 We next introduce the spaces of mesh functions associated with the
dual Delaunay/Voronoi grids. In an MFD scheme, the unknowns are
the components of $\mathbf u$ parallel to the
edges of the Delaunay triangulation and evaluated at the middle of
these edges. We orient each of the Delaunay edges by the unit vector
\[
{\mathbf e}_{ij}^D = {\mathbf e}_{ji}^D, \; i=1,\ldots,N_D,\;
j\in{\cal W}^V(i) = \{j\,\vert \, \partial V_{ij}\neq \emptyset,\,
j=1,\ldots,N_D\},
\]
which is directed from the node with the smaller index to the node
with  larger index, see Figure~\ref{de-a}.
This convention defines a function $\eta(i,j)$ for every edge
$(\mathbf x_{i}^D, \mathbf x_{j}^D)$,
\begin{equation}\label{eq:eta}
\eta(i,j) = \mathbf{n}_{i}^V\cdot \mathbf e_{ij}^D,
\end{equation}
where $\mathbf{n}_{i}^V$ is the unit normal vector outward to $ \partial V_i$.
We then denote by ${\mathbf H}_D$ the set of mesh functions
approximating the values of $\mathbf{u} \cdot {\mathbf e_{ij}^D}$ at
the mid-points of the edges connecting
$\mathbf x_i^D$ and $\mathbf x_j^D$; namely,
for all $i=1,\ldots N_D$ and $j$ such that $(\mathbf x_i^D, \mathbf
x_j^D)$ is a Delaunay edge we set
\[
u_{ij}^D\approx {\mathbf u}\cdot {\mathbf e_{ij}^D}(x_{ij}^D), \quad
{\mathbf x}_{ij}^D = \frac{1}{2}({\mathbf x}_{i}^D+{\mathbf x}_{j}^D).
\]

\begin{figure}[htb] \centering
\subfloat[]
{
\includegraphics*[scale = 0.4]{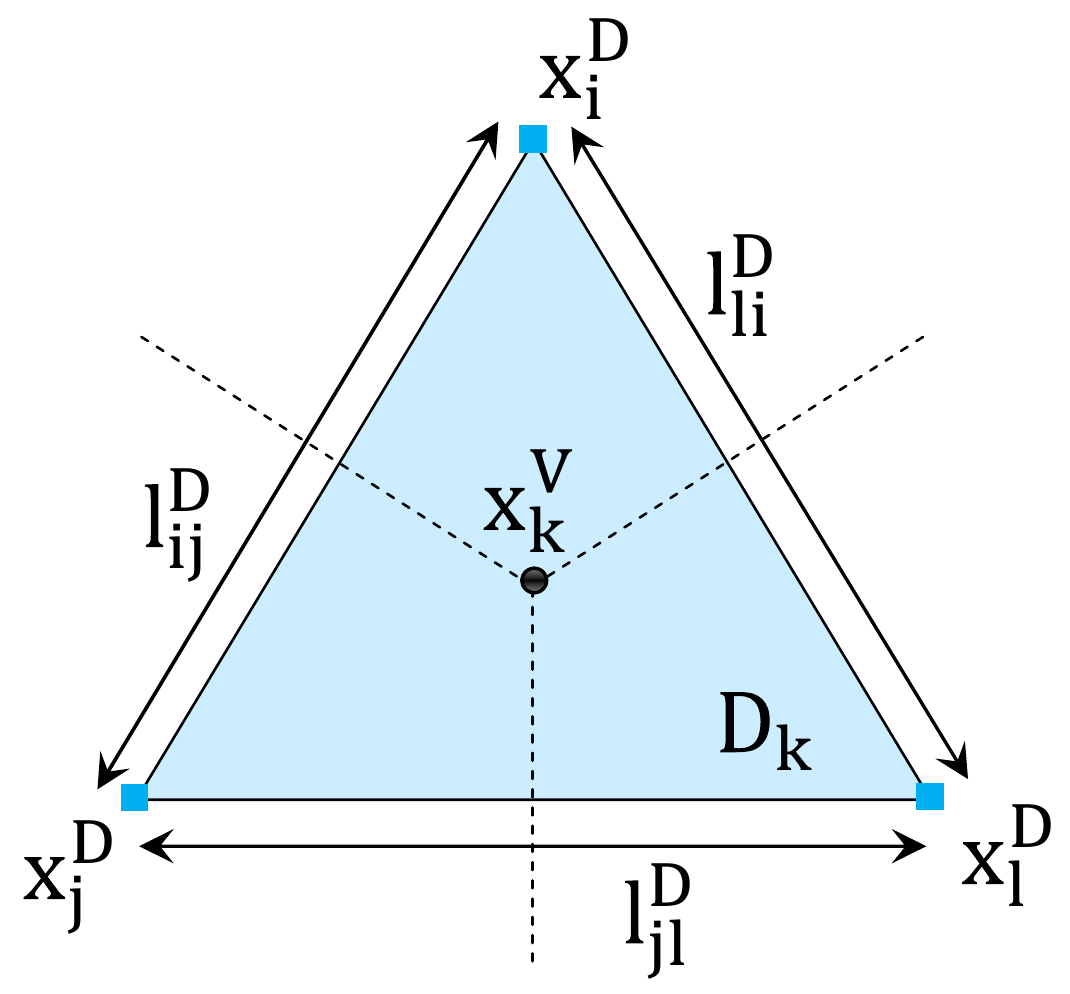}
\label{nr-a}
}
\subfloat[]
{
 \includegraphics*[scale =0.4]{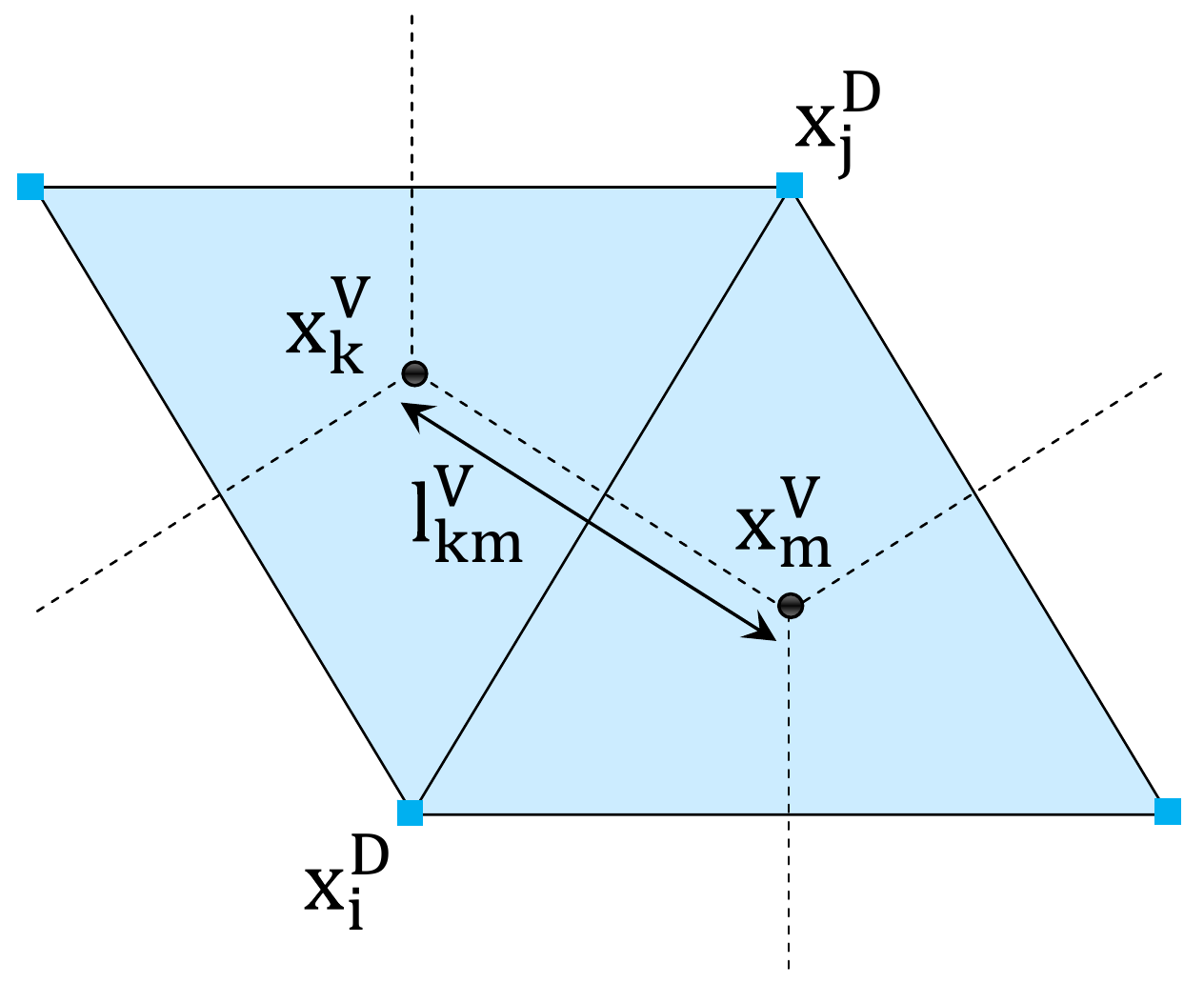}
\label{nr-b}
}
\caption{Notation used in the definition of $\rot_h$ and
$\curl_h$ operators}
\label{notation_rotor}
\end{figure}

To complete the MFD approximation of
equations~\eqref{rot_rot_problem}-\eqref{grad_div_problem} we now
introduce the spaces of mesh functions associated with the scalar
quantities $\rot \mathbf u$ and $\ddiv \mathbf u$.  We modify a little
bit of the definitions given in~\cite{Mimetic_Vabishchevich} to serve
better our purposes, although essentially we do not change anything
quantitatively.  With the vertices of the Delaunay (resp. Voronoi)
grid we associate the space of piece-wise constant functions,
which are constants on the polygons of the Voronoi (resp. Delaunay)
grid. We set
\begin{eqnarray}
&& H_D = \{u({\mathbf x})\, \vert \, u({\mathbf x})=u_i^D,\;\mbox{for
  all}\; x\in V_i, i=1,\ldots,N_D\},\\ \label{space-scalar-D}
&& H_V = \{u({\mathbf x})\, \vert \, u({\mathbf x})=u_k^V,\;\mbox{for  all}\; x\in D_k, k=1,\ldots,N_V\}. \label{space-scalar-V}
\end{eqnarray}
In short, the functions in $H_D$ are constants on Voronoi cells and the
functions in $H_V$ are constants on Delaunay cells.
We then define the discrete divergence operator as following:
\begin{equation}\label{eq:divh}
(\ddiv_h\, {\mathbf u})_i^D :=
\frac{1}{\meas(V_i)}\sum_{j: V_{ij}\in\partial V_i}u_{ij}^D\meas(V_{ij}),
\end{equation}
where $\meas(V_i)$ is the area (volume in 3D) of the Voronoi polygon
$V_i$, and, $\meas(V_{ij})$, is the length (area in 3D) of the Voronoi
edge (face in 3D) which is dual (perpendicular) to the Delaunay edge
$(\mathbf{x}^D_i,\mathbf{x}_j^D)$. The relation~\eqref{eq:divh} is
clearly an analogue of the divergence Theorem on $V_i$, namely,
\begin{equation*}
\frac{1}{\meas(V_i)}\int_{V_i} \ddiv {\mathbf u} =
\frac{1}{\meas(V_i)}\int_{\partial V_i} \mathbf{u}\cdot\mathbf{n}.
\end{equation*}
In a similar fashion we define the discrete  operators
$\grad_h: H_D \rightarrow \mathbf{H}_D$, $\rot_h: \mathbf{H}_V
\rightarrow H_D$, and
$\curl_h: H_V \rightarrow \mathbf{H}_D$
\begin{eqnarray}
&& (\grad_h\, u)_{ij}^D:=(\grad_h\, u)({\mathbf x}_{ij}^D)\cdot {\mathbf e}_{ij}^D = \eta(i,j) \frac{u_j^D-u_i^D}{l_{ij}^D},\label{eq:gradh}\\
&&(\rot_h \,{\mathbf u})^V_k = \frac{\eta(i,j)\,u_{ij}^D\,l_{ij}^D+\eta(j,l)\,u_{jl}^D\,l_{jl}^D+\eta(l,i)\,u_{li}^D\,l_{li}^D}{\hbox{meas}(D_k)},\label{eq:roth}\\
&&(\curl_h\, u)_{ij}^D = \eta(k,m)\frac{u_k^V - u_m^V}{l_{km}^V}, \label{eq:curlh}
\end{eqnarray}
where, in 2D, $\meas(D_k)$ is the area of the triangle with vertices
${\mathbf x}_i^D$, ${\mathbf x}_j^D$ and ${\mathbf x}_l^D$;
$l_{ij}^D$, $l_{jl}^D$ and $l_{li}^D$ are the lengths of its edges,
(for example, $l_{ij}^D=|x_i^D-x_j^D|$); and
$l_{km}^V  = |{\mathbf x}_k^V-{\mathbf
  x}_m^V|$. We refer to Figure~\ref{nr-a}-\ref{nr-b} for clarifying this notation.

\begin{figure}[htb] \centering
\centering
\subfloat[Stencil for $\curl_h\rot_h$]
{
\includegraphics*[width=0.48\textwidth]{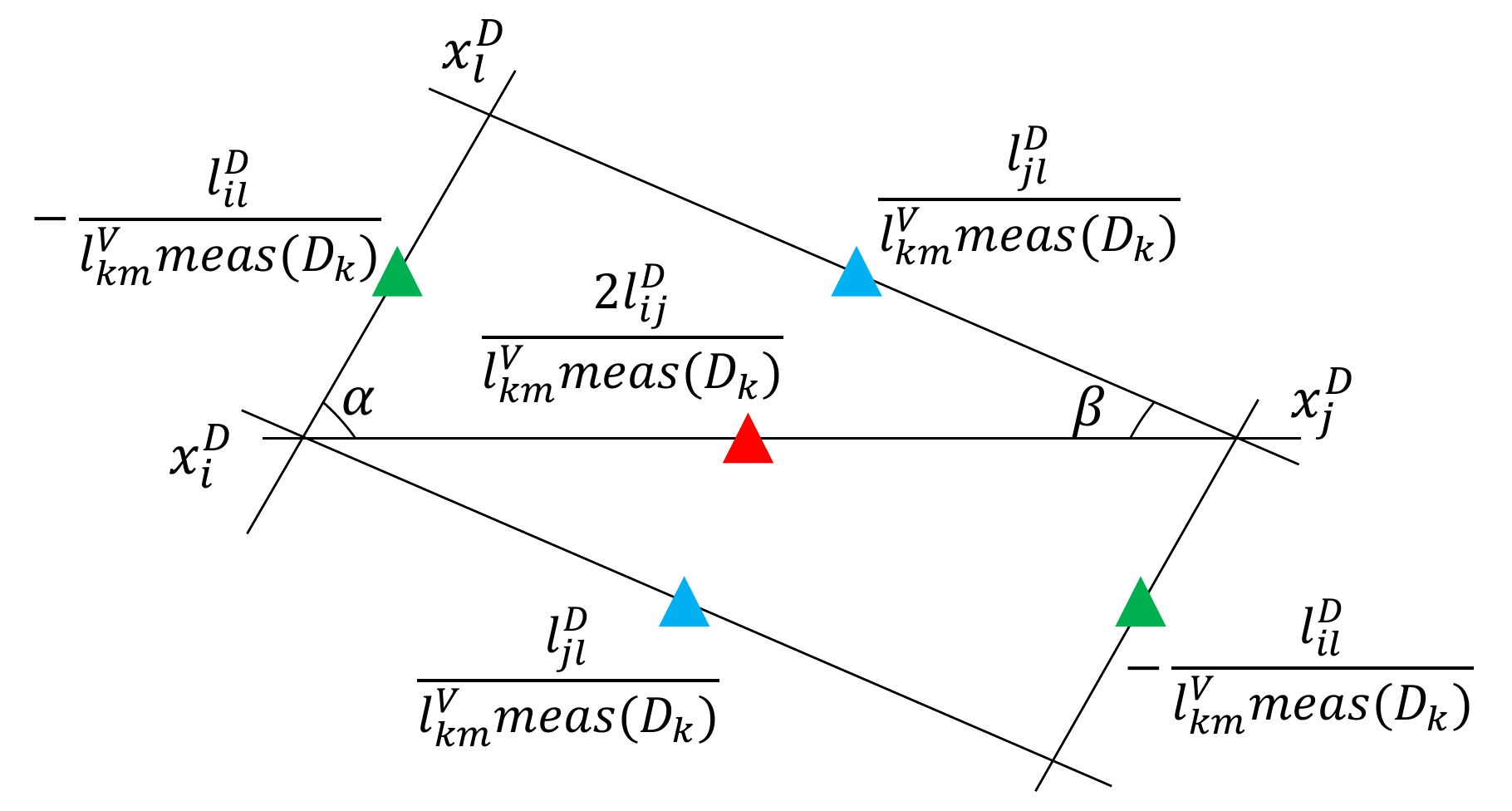}
\label{modelp-a}
}
\subfloat[Stencil for $(-\grad_h \, {\ddiv}_h)$]
{
\includegraphics*[width=0.48\textwidth]{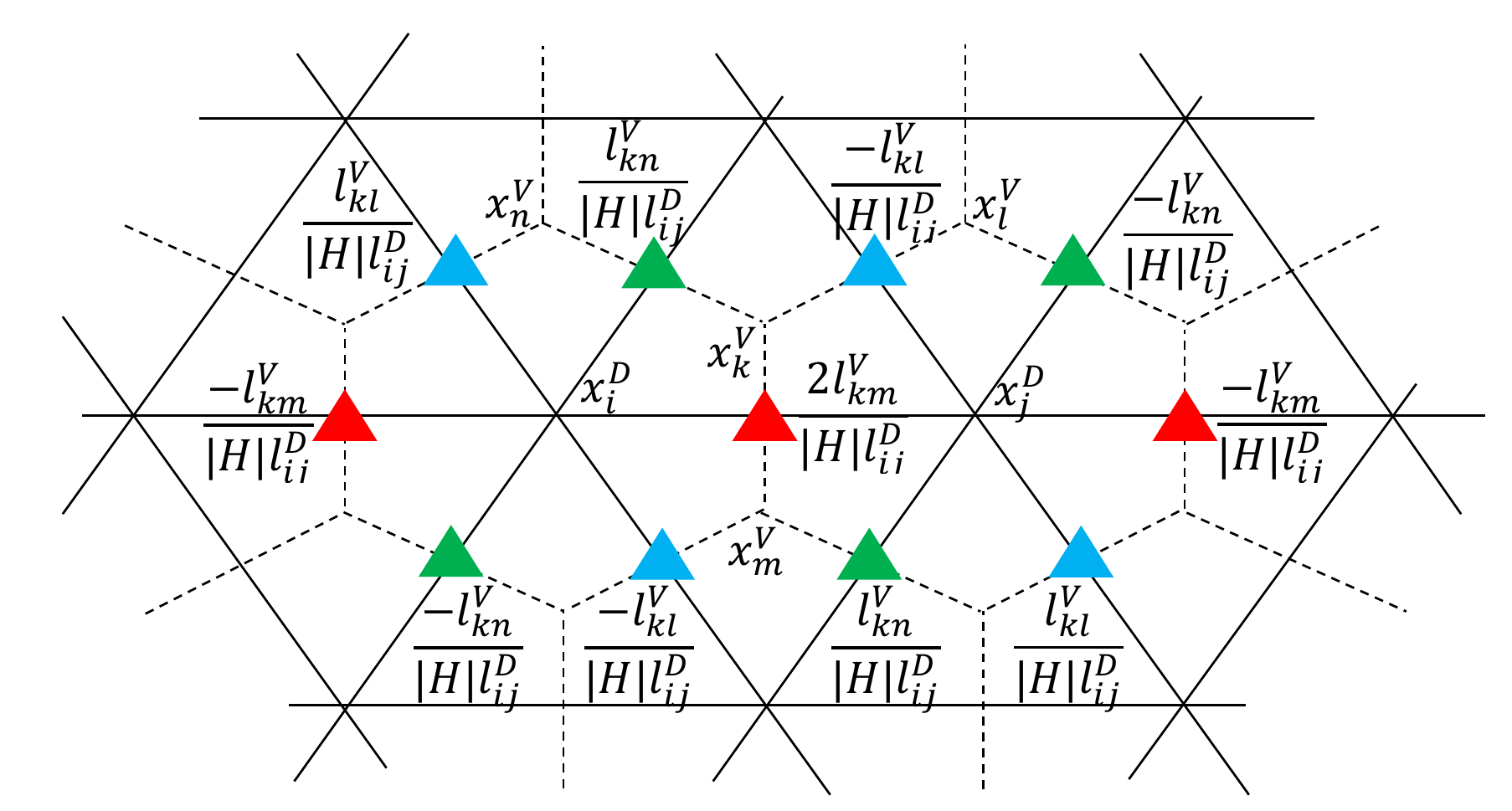}
\label{modelp-b}
}
\caption{Stencils corresponding to a uniformly refined triangular
mesh ($|H| = \meas(V_i)$)}
\label{modelproblems_general}
\end{figure}

Finally, The MFD stencils corresponding to the model
problems~\eqref{rot_rot_problem}--\eqref{grad_div_problem} are shown
in~Figure~\ref{modelp-a}-\ref{modelp-b} for a uniformly refined triangular grid.
Note that the stencils match exactly the ones given
in~\cite{Mimetic_Vabishchevich}. The modifications for the 3D variants of the operators above can also be found in~\cite{Mimetic_Vabishchevich}

\section{Equivalence between mimetic finite differences and finite element methods}\label{sec:FD_FEM}


In this section, for both model problems~\eqref{rot_rot_problem}
and~\eqref{grad_div_problem}, we are going to introduce suitable
FE methods to derive stencils on arbitrary structured
triangular grids that match those obtained by the MFD schemes. The FE methods that we consider are based on the variational
formulations~\eqref{rot_rot_problem-w} and~\eqref{grad_div_problem-w}. For
definitions of the corresponding Hilbert spaces and results on
existence and uniqueness of solutions to these model problems, we refer
to~\cite{Girault_Raviart}.

\subsection{Finite element discretization for~\eqref{rot_rot_problem-w}}
Next, we recall that the mesh functions are defined as approximations
to the tangential components of the solution on the Delaunay mesh. Therefore, if we would like to construct a FE discretization that matches the MFD method from the previous section, it is
reasonable to use lowest order $\Hcurl$-conforming
N\'ed\'elec elements~\cite{Nedelec1980,Nedelec1986} which have the $0$-th order moments on the edges of the Delaunay mesh of the tangential components of $\mathbf u$ as
degrees of freedom.

In order to approximate the variational
problem~\eqref{rot_rot_problem-w} by the lowest-order
N\'{e}d\'{e}lec's edge elements, we consider vector valued functions
whose restrictions on every Delaunay triangle $D_k$ are linear in each component and have tangential
components that are continuous across the element boundaries. Namely,
we define the N\'{e}d\'{e}lec's FE space as following
 \begin{equation}\label{Vh} {\mathbf V}^N_h = \{{\mathbf v}_h\in \Hcurl \, | \, {\mathbf v}_h|_{D_k} =
   \left[\begin{array}{c}a_1 \\ a_2 \end{array}\right]+b
   \left[\begin{array}{c}y \\ -x \end{array}\right],\, k=1,2, \dots, N_V \}. \end{equation}
The FE approximation of~\eqref{rot_rot_problem-w} is: Find $\mathbf{u}_h\in {\mathbf V}^N_h$ such
 that
\begin{equation}\label{discrete_weak_form}
   (\rot {\mathbf u}_h,\rot {\mathbf v}_h)
   +\kappa ( {\mathbf u}_h, {\mathbf v}_h)
   = ({\mathbf f},{\mathbf v}_h), \quad \forall
   \,{\mathbf v_h}\in {\mathbf V}^N_h.
\end{equation}

 As is well known, the degrees of freedom (functionals which uniquely
 determine the elements in $\mathbf{V}^N_h$) are chosen to ensure
 tangential continuity between elements and in the lowest order case,
 the degrees of freedom are the $0$-th order moments of the tangential component
 on each edge, i.e.
 $u_{ij} =\int_{x_{i}^D}^{x_{j}^D}{\mathbf u}_h\cdot {\mathbf
   e}_{ij}^D,$
 with ${\mathbf e}_{ij}^D$ defined as in Section~\ref{sec:mimeticFD}.

\begin{figure}[!htb]
\centering
\subfloat[N\'ed\'elec finite elements]
{
\includegraphics*[width=0.48\textwidth]{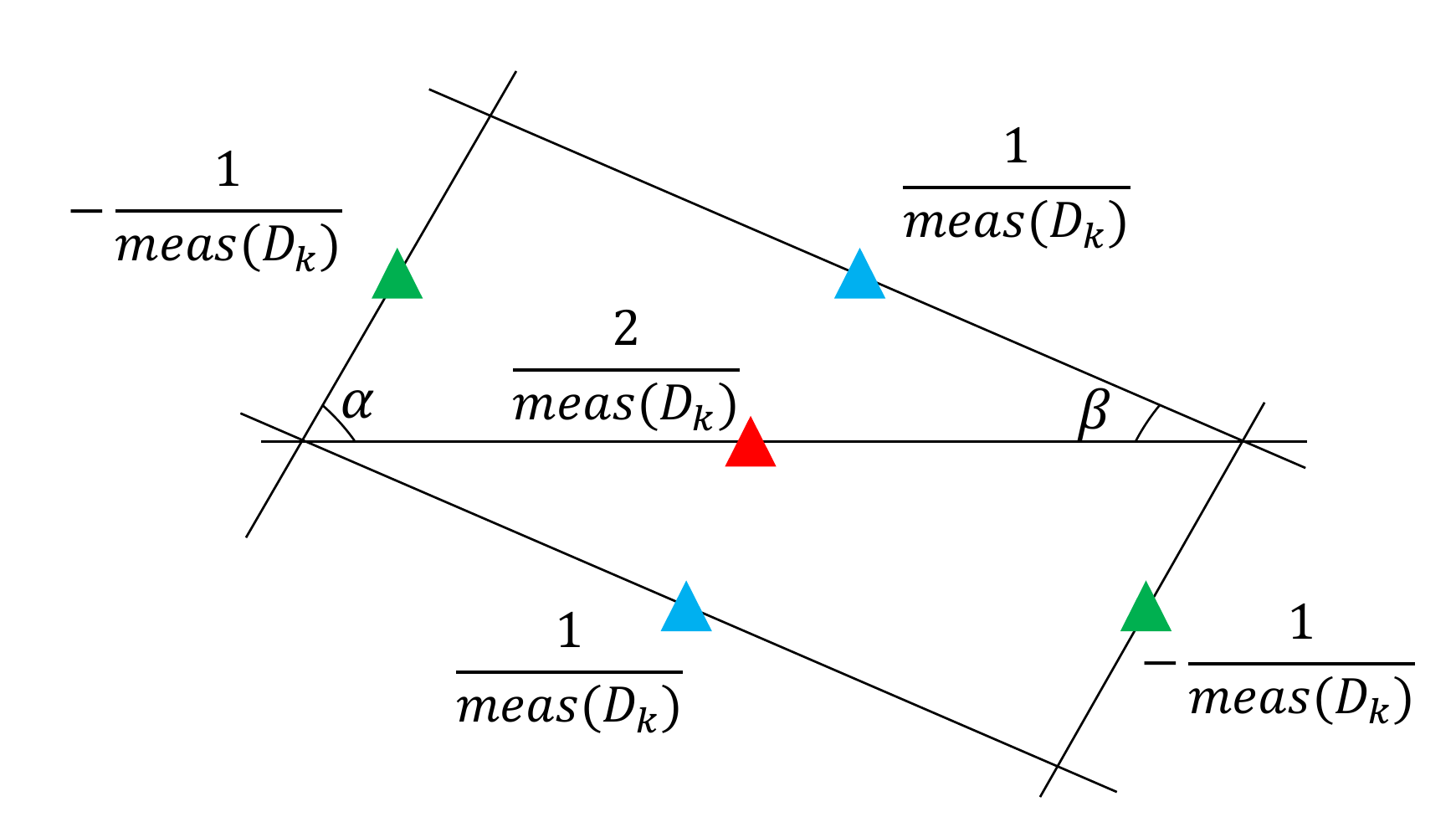}
\label{g-stencil-a}
}
\subfloat[Mimetic finite  differences]
{
\includegraphics*[width=0.48\textwidth]{eq_rot_rot_general_2}
\label{g-stencil-b}
}
\caption{Stencils for a general triangulation.}
\label{general_stencil_2}
\end{figure}

 The bases dual to these degrees of freedom have one basis function
 per Delaunay edge $(\mathbf x_i^D, \mathbf x_j^D)$,
 $\varphi_{ij} = \frac{1}{2}(\lambda_i \nabla \lambda_j - \lambda_j
 \nabla \lambda_i)$,
 where $\lambda_i$ and $\lambda_j$ are the barycentric coordinates of the
 Delaunay grid.  In the standard fashion, the solution of
 problem~\eqref{discrete_weak_form} is written as
${\mathbf u}_h = \displaystyle\sum_{(i,j)}u_{ij}\varphi_{ij}$,
and the vector of coefficients $U^N = (u_{ij})$ is a solution to the
linear system of equations $A^N U^N = b^N$. Here,
the stiffness matrix $A^N$ has elements given by
\begin{equation}
(A^N)_{(i_2,j_2)(i_1,j_1)}=\sum_{k=1}^{N_V}\int_{D_k}\left( \rot\varphi_{i_1 j_1} \,
\hbox{rot} \varphi_{i_2 j_2} + \kappa \, \varphi_{i_1 j_1}
\varphi_{i_2 j_2} \right)\mathrm{d}\mathbf{x}.
\end{equation}
We can now compare the stencil (a row in $A^N$) corresponding to the
FE discretization and the MFD discretization derived in~Section~\ref{sec:mimeticFD}.
The result can be seen in Figure~\ref{g-stencil-a}--\ref{g-stencil-b} and there is obviously no match.
In order to see the relation between the MFD and the FE stencils, the
key is that the MFD degree of freedom $u_{ij}^D$ is basically
a scaled N\'{e}d\'{e}lec FE degree of freedom $u_{ij}$, and since
the degrees of freedom are dual to the basis functions, we need to
scale appropriately the basis in order to have the same entries in
the FE stencil.  More precisely, if
we use the midpoint quadrature rule on every edge (which is exact for
functions in $\mathbf{V}^N_h$) we have
\begin{equation}\label{u_via_nedelec}
{\mathbf u}_h({\mathbf x}) =
  \displaystyle\sum_{(i,j)}\left(\int_{x_{i}^D}^{x_{j}^D}{\mathbf
      u}_h\cdot {\mathbf e}_{ij}^D\right)\varphi_{ij}({\mathbf x}) =
  \displaystyle\sum_{(i,j)}({\mathbf u}_h\cdot {\mathbf
    e}_{ij}^D)(x_{ij}^D)\, l_{ij}^D\, \varphi_{ij}({\mathbf
    x}).
\end{equation}
Hence, an appropriate re-scaling of the basis is
$\varphi_{ij}^{s} := l_{ij}^D\, \varphi_{ij}$. Using the re-scaled
basis functions, the corresponding stencil is still not in agreement
with the MFD but the resulting rows are proportional to each other. We
easily remedy this by scaling also the test functions
$\widetilde{\varphi}_{ij}=\displaystyle\frac{1}{l_{km}^{V}}\,\varphi_{ij}$. Summarizing,
the MFD matrix $A^{FD}$ and the FE matrix $A^N$ satisfy
\begin{equation}\label{relation_rot_rot}
  A^{FD} = D_1\, A^N \, D_2, \quad \hbox{where}
  \left\{\begin{array}{l} D_1 = \diag((l_{km}^V)^{-1}) \\ D_2 =
      \diag(l_{ij}^D) \end{array} \right.
\end{equation}

For the right hand side $\mathbf{f}$, in order to be consistent with the FE method, we use the following approximation
\[
f_{ij}^D = \frac{1}{l_{km}^V} \int_{\Omega} \mathbf{f} \varphi_{ij} \mathrm{d} \mathbf{x}.
\]
Then we have $b^{FD} = D_1 b^N$ and, finally, obtain the following relation
between the MFD solution $U^{FD} = (u_{ij}^D)$ and the N\'ed\'elec FE solution given by $U^{N}$:
\begin{equation}\label{all-ok}
A^{FD} U^{FD} = b^{FD}
\Longrightarrow (D_1 A^N D_2) (U^{FD}) = D_1 b^N, \quad
U^{FD} = D_2^{-1} U^{N}.
\end{equation}
Clearly, there is a function in~$\mathbf{V}^N_h$ corresponding to the MFD solution which we can define as following (in the notation of Section~\ref{sec:mimeticFD})
\begin{equation*}
  \mathbf{u}_h^{FD}(\mathbf{x}) =
  \sum_{(i,j)} u^{D}_{ij}
  \varphi_{ij}^{s}(\mathbf{x}).
\end{equation*}
Using the fact that $U^{FD} = D_2^{-1} U^{N}$ and $\varphi_{ij}^{s}=l_{ij}^{D}\,\varphi_{ij}$ we obtain that $\mathbf{u}_h^{FD}=\mathbf{u}_h$.
As a consequence, from the standard error analysis for the N\'ed\'elec FE methods
for sufficiently regular $\Omega$ we automatically have  an error estimate and stability for the MFD discretization, namely,
\begin{equation}\label{error-rot-rot}
\|\mathbf{u} - \mathbf{u}_h^{FD} \|_{\mathrm{rot}} \leq C h \|
\mathbf{f} \|_{L^2(\Omega)},
\end{equation}
where
$\| \mathbf{v} \|_{\rot} := \sqrt{(\rot\mathbf{v},
  \rot\mathbf{v})+\kappa(\mathbf{v},\mathbf{v})}$
and $\| \cdot\|_{L^2(\Omega)} $ is the standard $L^2$-norm on
$\Omega$.

We can also use other approximations of $\mathbf{f}$ in the MFD schemes.  This implies to use an approximation $\widetilde{\mathbf{f}}$ in~\eqref{discrete_weak_form} instead of $\mathbf{f}$.  It is reasonable to assume that $\| \mathbf{f} - \widetilde{\mathbf{f}} \|_{L^2(\Omega)} \leq C h $, and, therefore, we still have the error estimate~\eqref{error-rot-rot} by the standard perturbation argument and triangular inequality.

\begin{remark}\label{remxxx}
 Although the matrix $A^{FD} = D_1A^N D_2$ could
  be non-symmetric in the standard Euclidean inner product, it is symmetric in the inner product defined by
  $(D_1^{-1} D_2)$. This is a crucial observation which plays an important
  role in the design of efficient solvers for the MFD discretizations
  of problem~\eqref{rot_rot_problem}.
\end{remark}

\subsection{Finite element methods for~\eqref{grad_div_problem}}
\begin{figure}[!htb]
\centering
\includegraphics[scale =0.56]{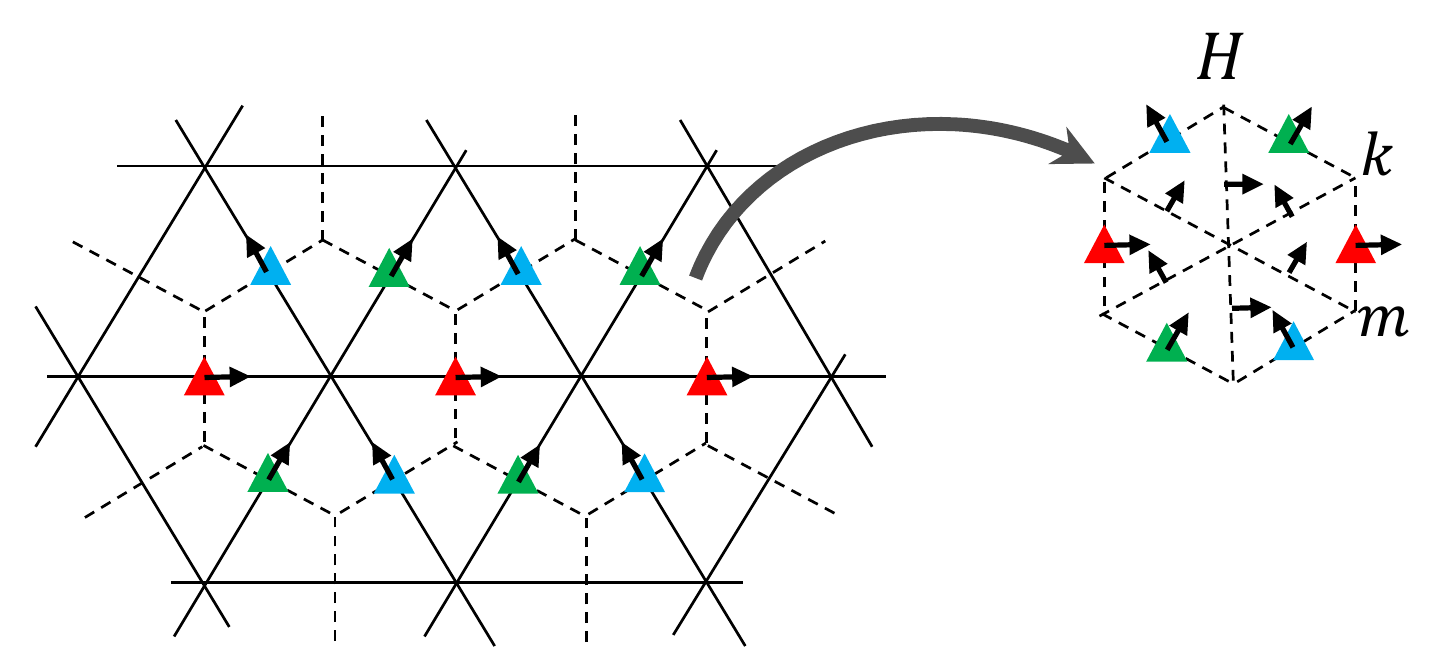}
\caption{Splitting of an hexagon into six
triangles to compute the Raviart-Thomas basis function on the
hexagon.\label{RT_basis_hexagon} }
\end{figure}

The situation with the model problem~\eqref{grad_div_problem} is a bit
more involved. To obtain the FE discretization that matches the MFD
discretization of~\eqref{grad_div_problem} we borrow some ideas
from~\cite{Kuznetsov2003,Kuznetsov2004,Pasciak_Vassilevski} and
construct an FE space on the polyhedral grid formed by the
Voronoi cells. Let us recall that we have $N_D$ Voronoi cells. We consider the space of functions whose divergence is constant on these cells.
 \begin{equation}\label{VhRT} {\mathbf V}^{\text{RT}}_h =
   \{{\mathbf v}_h\in \Hdiv \, | \, {\mathbf v}_h|_{V_i} \in \Hdiv(V_i), \; \ddiv \mathbf v_h = \text{const}.
   \}.
\end{equation}
While the particular behavior of the functions inside the Voronoi cell
does not matter for our considerations that follow, just to fix the
space, we assume that the elements of $\mathbf{V}^{\text{RT}}_h$ are
in the standard Raviart-Thomas space~\cite{Raviart_Thomas_original_paper,Nedelec1986}
on a sub-triangulation of every Voronoi cell as shown in
Figure~\ref{RT_basis_hexagon} with fixed values of the
Raviart-Thomas degrees of freedom on $\partial V_i$ for
$i=1,\ldots, N_D$.  With this choice, the finite element approximation
of~\eqref{grad_div_problem-w} is: Find
$\mathbf{u}_h\in {\mathbf V}^{\text{RT}}_h$ such that
\begin{equation}\label{discrete_weak_form-1}
   (\hbox{div}\, {\mathbf u}_h,\hbox{div}\, {\mathbf v}_h)
   +\kappa ( {\mathbf u}_h, {\mathbf v}_h)
   = ({\mathbf f},{\mathbf v}_h), \quad \forall
   \,{\mathbf v_h}\in {\mathbf V}^{\text{RT}}_h.
\end{equation}
 There is one degree of freedom associated with
each Voronoi edge $(\mathbf x_{k}^V, \mathbf x_m^V)$. The
corresponding basis functions, dual to these degrees of freedom, are
defined using a sub-triangulation of the hexagon as shown in
Figure~\ref{RT_basis_hexagon} (although this is not necessary,
see~\cite{Kuznetsov2004}), and then solving an auxiliary finite
element problem in each hexagon.
We denote by $RT(H)$ the standard Raviart-Thomas space corresponding to
a sub-triangulation of a Voronoi cell $H$. If we are on a uniformly
refined grid, this is a hexagon, which explains why we denoted a
generic cell with $H$, but in general it is just a convex cell. We
restrict our considerations to a hexagon $H$, since the other cases
are similar.  The basis function on $H$ corresponding to the edge
$(\mathbf x_{k}^V, \mathbf x_m^V)\subset \partial H$ is determined by
solving the following constrained minimization problem: Find
$\varphi_{km}\in RT(H) \; \; (dim\, RT (H) = 12)$ such that
\renewcommand{\arraystretch}{1.5}
\[
\left\{
\begin{array}{l}
\|\varphi_{km}\|^2_{*} \rightarrow \min{}, \; \|\varphi_{km}\|_* \simeq \|\varphi_{km}\|_{L^2},\;\;
\ddiv\varphi_{km} = \pm \frac{1}{\meas(H)}, \\
                  \displaystyle\int_{\mathbf x_{k}^V}^{\mathbf x_m^V}
\varphi_{km}\cdot {\mathbf n}_{jl} = \delta_{(km),(jl)},\; \forall jl\in\partial H,
\end{array}
\right.
\]
where $\|\cdot\|_*$ can be the $L^2(H)$-norm
or any equivalent norm on the space $RT(H)$.

\renewcommand{\arraystretch}{1}

As in the FE discretization of~\eqref{rot_rot_problem}, which we
considered above, to match the MFD discretization, we need to re-scale
both the test and trial functions. We set
$\varphi_{km}^{s} = l_{km}^{V}\varphi_{km}$, and take as new test
functions $\widetilde{\varphi}_{km} = \frac{1}{l_{ij}^D}\varphi_{km}$.
We then obtain that \begin{equation}\label{matrices_grad_div}
  A^{FD} = D_1\, A^{RT} \, D_2,\qquad \hbox{where} \;\left\{
\begin{array}{l} D_1 =
                                          \diag((l_{ij}^D)^{-1}) \\
 D_2 = \diag(l_{km}^V) \end{array}
                                      \right.
\end{equation}
Similarly, properly defining the approximation of
$\mathbf{f}$, we have
\begin{equation*} A^{RT} U^{RT} = b^{RT} \ \text{and} \ A^{FD}
U^{FD} = b^{FD}
\end{equation*}
with $b^{FD} = D_1 b^{RT}$.
Together with~\eqref{matrices_grad_div}, we have $U^{FD} = D_2^{-1} U^{RT}$.
Therefore, the discrete MFD solution corresponds to the following function in $\mathbf{V}_h^{\text{RT}}$
\begin{equation*}
\mathbf{u}_h^{FD}(\mathbf{x}) = \sum_{(k,m)} u^{D}_{km}
\varphi_{km}^{s}(\mathbf{x}), \end{equation*}
which also satisfies $\mathbf{u}_h^{FD}=\mathbf{u}_h$, namely the MFD solution is also the FE solution.
Applying then some standard arguments for a sufficient regular domain
$\Omega$, we have the following error estimate
\begin{equation} \label{error-grad-div}
\| \mathbf{u} - \mathbf{u}_h^{FD} \|_{\mathrm{div}} \leq C h \| \mathbf{f}
\|_{L^2(\Omega)}, \end{equation}
where $\| \mathbf{v}
\|_{\mathrm{div}} := \sqrt{(\mathrm{div} \ \mathbf{v}, \mathrm{div} \
\mathbf{v}) + \kappa (\mathbf{v}, \mathbf{v})}$.

We can also use other approximation $\widetilde{\mathbf{f}}$ of $\mathbf{f}$ and, similarly, obtain the error estimate~\eqref{error-grad-div} by standard perturbation argument and triangular inequality.  Moreover, we note that a remark analogous to Remark~\ref{remxxx} is necessary here as well.

\section{Multigrid solvers for mimetic finite differences}\label{sec:FEM_mg}

Our aim is to find efficient multigrid (MG) methods for the MFD discretization of the vector problems. In this section, we focus on
problem~\eqref{rot_rot_problem}.
A geometric multigrid (GMG) method for the MFD discretizations on the Voronoi cells (hexagonal grids) is
a topic of our ongoing research and will be reported in our future work.  We point out, however, the FE framework that we established
in the previous section provides the necessary conditions for applying
efficient methods using irregular coarsening strategies and algebraic
multigrid together with auxiliary space methods
(see~\cite{Lashuk_Vassilevski2012,Hiptmair2007}). Such methods,
however, are not suitable for local Fourier analysis.  Since one of
our goals is to perform the LFA, we focus on the GMG method
for problem~\eqref{rot_rot_problem} here.

\begin{figure}
\includegraphics[scale = 0.6]{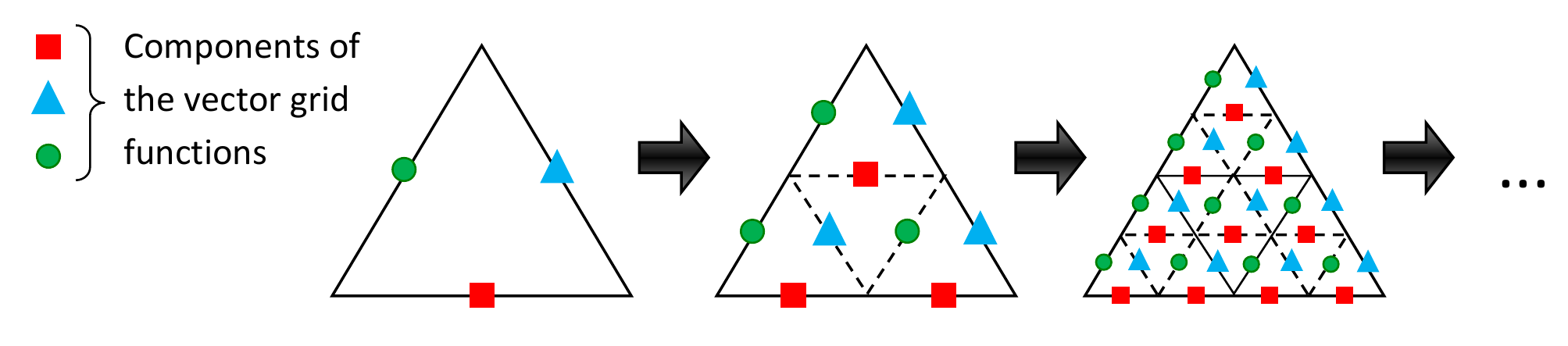}
\caption{Hierarchy of grids to perform the geometric multigrid method,
and location of the unknowns}
\label{hierarchy}
\end{figure}

We are interested in applying a GMG method on
triangular grids generated by regular refinement. In this way, we
naturally obtain a hierarchy of grids, as shown in~Figure~\ref{hierarchy}.
Next we describe
the components for the MG algorithm, i.e., smoother and inter-grid
transfer operators.

\subsection{Multigrid components}\label{sec:mg_components}
\begin{figure}
\centering
\includegraphics[width=0.3\textwidth]{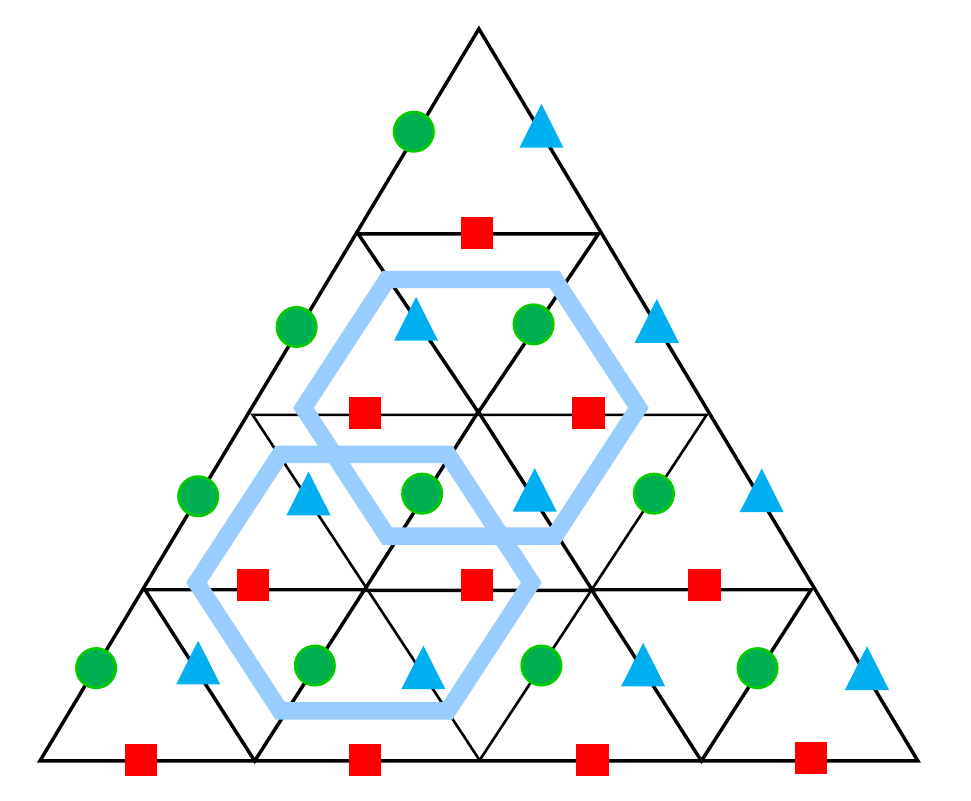}
\caption{Unknowns simultaneously updated in the overlapping block
smoother, and overlapping of the blocks.}
\label{vanka_smoother}
\end{figure}
\paragraph{Smoother}
We use a multiplicative Schwarz smoother proposed in~\cite{Arnold_et_al} as the relaxation (smoother) in the standard $V$- and $W$-cycle.  This
relaxation simultaneously updates all the unknowns around a vertex of
the Delaunay grid as shown in~Figure~\ref{vanka_smoother}. Overlapping of the
unknowns requires non-standard tools in the local Fourier analysis for this smoother.

\begin{figure}[htb] \centering
\subfloat [Canonical N\'ed\'elec prolongation matrix.]
{
\includegraphics[scale = 0.65]{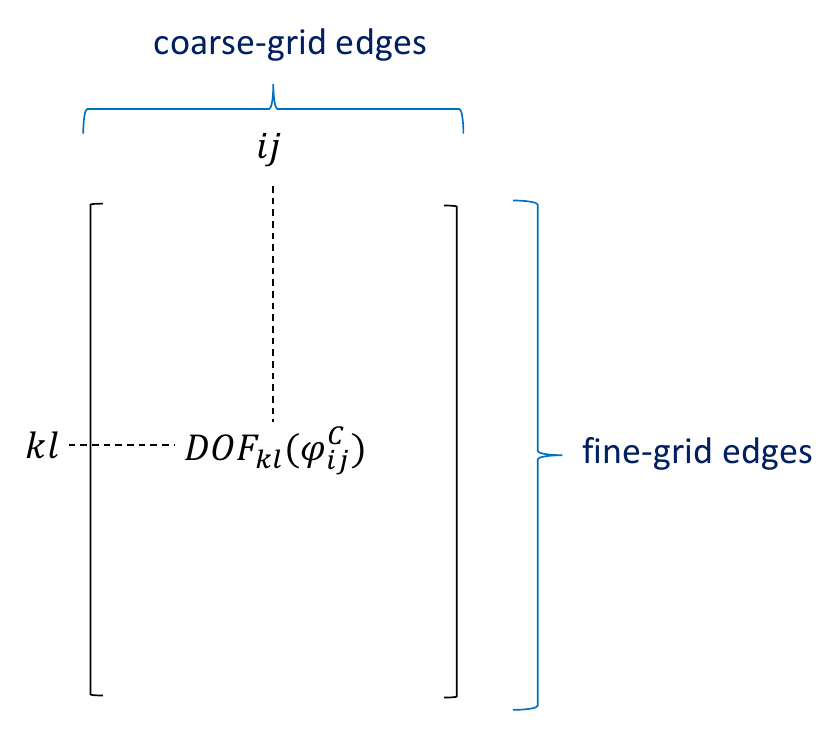}
\label{N-a}
}\hspace{2em}
\subfloat[Weights for the canonical N\'ed\'elec prolongation operator.]
{
\includegraphics[scale = 0.55]{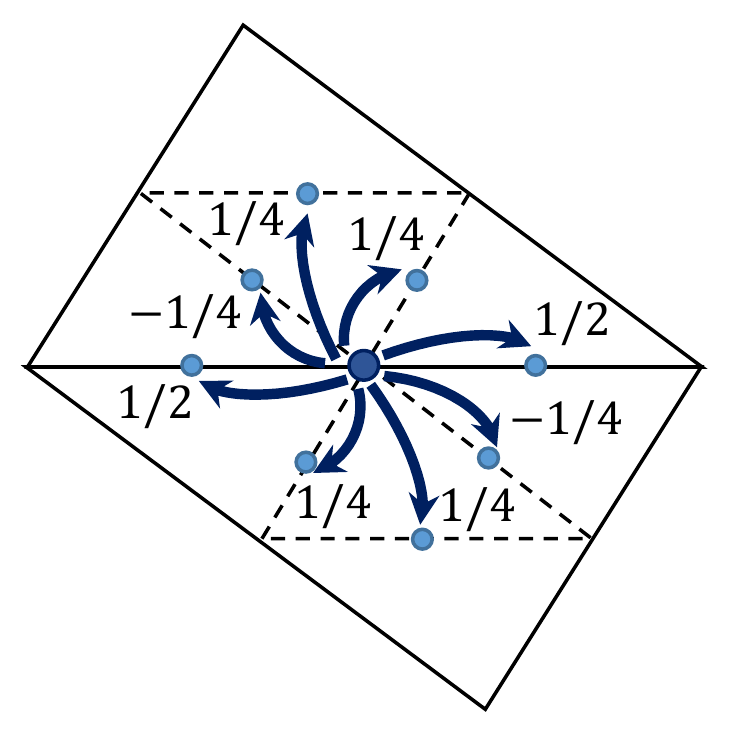}
\label{N-b}
}
\caption{}
\label{Nedelec_prolongation}
\end{figure}

\begin{figure}
\centering
\includegraphics[width=0.48\textwidth]{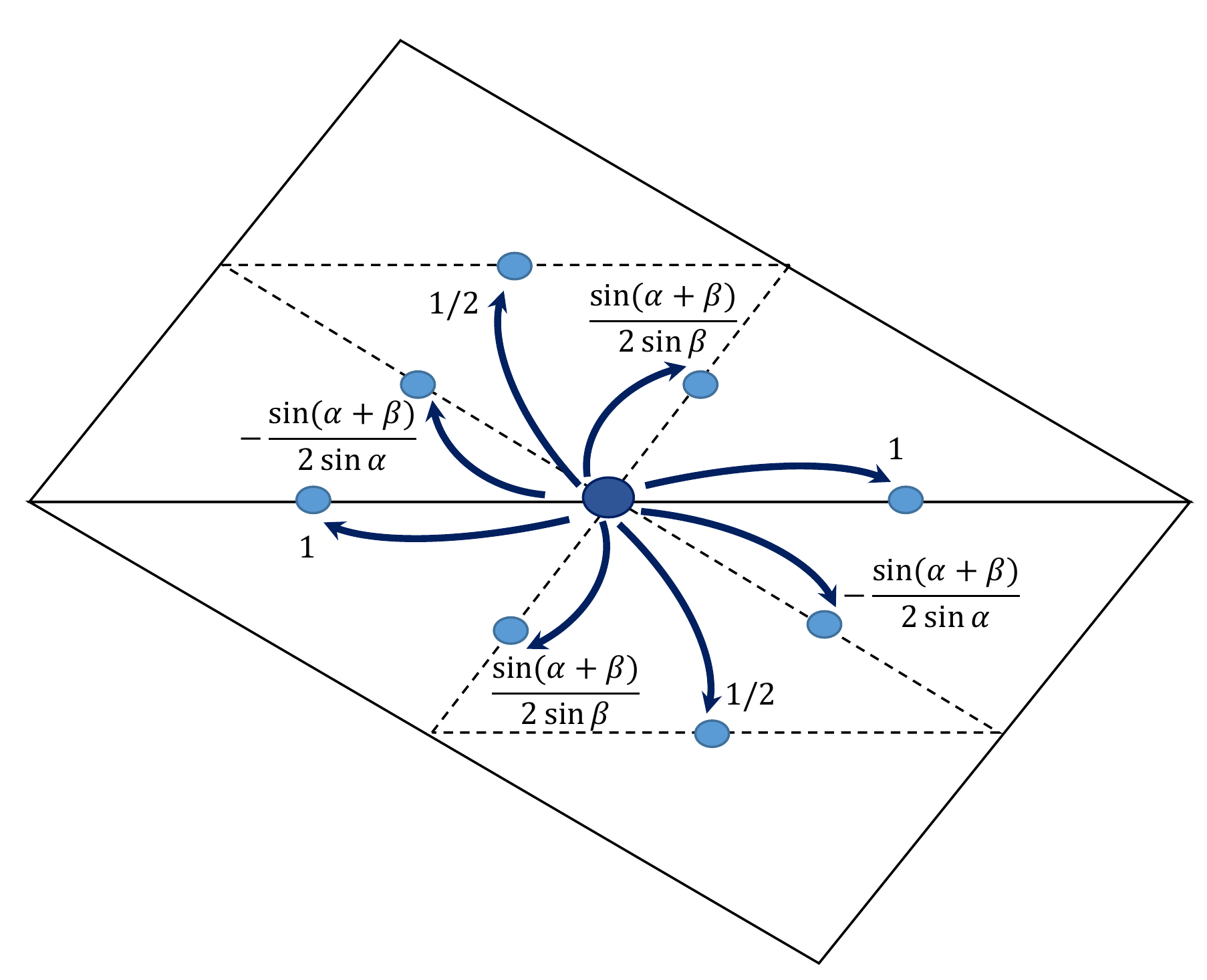}
\caption{Weights for the canonical prolongation associated with the
``modified'' N\'ed\'elec finite element scheme.}
\label{canonical_prolongation}
\end{figure}

\paragraph{Inter-grid transfer operators} For the transfer of
information between two consecutive grids of the hierarchy, canonical
inter-grid transfer operators based on the FE framework are
constructed.  We remind that in order to obtain the N\'ed\'elec
canonical prolongation, we need to write the coarse-grid basis
functions as a linear combination of the fine-grid basis functions,
that is,
$$ \quad \varphi_{ij}^C = \displaystyle \sum_{kl\in{\mathbf K}} DOF_{kl}(\varphi_{ij}^C)\varphi_{kl}^F,$$
where ${\mathbf K}$ denotes the set of fine-grid edges inside the
support of $\varphi_{ij}^C$.
In this way, the coefficients in the
linear combination are the entries of the prolongation matrix $P^N$, i.e. $(P^N)_{(k,l),(i,j)} = DOF_{kl}(\varphi_{ij}^C)$,
see Figure~\ref{N-a}. Moreover, the weights in $P^N$ are displayed in Figure~\ref{N-b}.

Recall that we needed to re-scale the standard  N\'ed\'elec basis in
order to obtain the equivalence between the MFD and FE methods.
We modify accordingly the prolongation operator and, from N\'ed\'elec canonical prolongation $P^N$,
we can write everything in terms of the re-scaled
fine-grid and coarse-grid basis functions to obtain that
the prolongation matrix is $P=D_{2,h}^{-1}P^N D_{2,H}$.
The resulting prolongation for a general
refined triangle with angles $\alpha$ and $\beta$ is given in
Figure~\ref{canonical_prolongation}. 
For the restriction, by similar re-scaling,
we choose
$$R = D_{1,H} (P^N)^T D_{1,h}^{-1}.$$
It is worth noting here that $R$ is the adjoint of $P$ in the inner
products induced by $(D_{1,h}^{-1}D_{2,h})$ and $(D_{1,H}D_{2,H}^{-1})$ on the
fine and coarse grids, respectively, and this agrees with the
considerations in Remark~\ref{remxxx}.

The choice of these inter-grid transfer operators is, in fact, crucial. We emphasize
the importance of the relation obtained in Section~\ref{sec:FD_FEM}
because without this, it is possible, but, by all means not easy to
design efficient GMG methods for these MFD discretizations.
Regarding the coarse-grid MFD operator, with this choice
of the inter-grid transfer operators, direct
discretization on the coarse-grid results in a ``Petrov-Galerkin''
coarse-grid operator. Namely, note that since $A_H^N
= (P^{N})^TA_h^N P^N$, we have
\begin{eqnarray*} A_H^{FD} &=& D_{1,H}A_H^N D_{2,H} = D_{1,H}(P^N)^T A_h^N
P^N D_{2,H} \\ &=& D_{1,H}(P^N)^T \overbrace{D_{1,h}^{-1} D_{1,h}}^{I}A_h^N
\overbrace{D_{2,h}D_{2,h}^{-1}}^{I} P^N D_{2,H} \\ &=&
\overbrace{(D_{1,H}(P^{N})^T D_{1,h}^{-1})}^{R} A_h^{FD}
\overbrace{(D_{2,h}^{-1} P^N D_{2,H})}^{P} \\ &=& RA_h^{FD}P.
\end{eqnarray*}

\subsection{Local Fourier analysis}\label{sec:LFA}

We now briefly use local Fourier analysis techniques to assess the
convergence of the resulting GMG method. The LFA (or
local \emph{mode}) analysis is introduced by A.~Brandt
in~\cite{Brandt77} and is a technique based on the Discrete Fourier
Transform.
To perform this analysis one slightly diverts from the boundary
value problem in
hand and considers periodic solutions on an infinite
regular grid. It is also necessary to have a discrete operator
defined with a constant coefficient stencil. The boundary conditions are
not taken into account, that is, we assume that the boundary effect is negligible. This, of course, is
not true in general, but in many practical situations the LFA gives
sharp estimates on the MG convergence rates.

In the framework of LFA, the current approximation to the solution and
the corresponding error can be represented by formal linear
combinations of discrete \emph{Fourier modes} forming the
\emph{discrete Fourier space}.
The LFA then identifies invariant subspaces in the Fourier space
and studies how multigrid components
act on these subspaces. 
A detailed explanation of all varieties of local Fourier analysis can be found
in~\cite{TOS01,Wie01} and on triangular grids, in~\cite{RodrigoBook}.

We now move on to describe the difficulties in performing LFA for the
GMG components defined above.
The LFA for our case is indeed nonstandard and we need to deal with several issues
described below.

\paragraph{Simplicial grids}
Local Fourier analysis has been traditionally performed for finite
difference discretizations on structured rectangular grids. This
analysis was extended to FE discretizations on general
structured triangular~\cite{LFA_tri, LFA_tri_mma} and
tetrahedral~\cite{LFA_tri_tetra} grids. The key fact for this
extension is to consider an expression of the Fourier transform in new
coordinate systems in space and frequency variables and introduce a
non-orthogonal unit basis of ${\mathbb R}^d$, chosen to fit the
geometry of the given simplicial mesh. The basis corresponding to the
frequencies space is taken as its reciprocal basis and with these
settings the LFA on simplicial grids is not very different from the
LFA on the standard rectangular grids.

\paragraph{Edge-based unknowns}
As we saw, the unknowns in MFD discretizations of
problem~\eqref{rot_rot_problem} are located at different types of
grid-points, and therefore the stencils (the rows of the matrix
$A^{FD}$) involve not one, but several different stencils.  The key is
to split the infinite grid into several different subgrids in such a
way that all nodes belonging to a subgrid have the same stencil,
and to define suitable grid-functions playing the role of the Fourier
modes for such edge-based discretizations. We refer the reader
to~\cite{RodrigoWeizmann} for a detailed description of such
analysis.

\paragraph{Overlapping Schwarz smoothers}
Overlapping block smoothers require a special LFA strategy.  Classical approaches fail for this class of smoothers, because an
overlapping smoother updates some variables more than once, due to the
overlapping. The main difficulty is that in addition to the initial and final
errors, some intermediate errors appear, and this has to be taken into account in the analysis.
To our knowledge, there are only few
papers dealing with LFA for overlapping smoothers, and all of them for
discretizations on rectangular grids (see~\cite{Oost-2010},
\cite{molenaar-91} and \cite{Siv-91}). In~\cite{RodrigoBook} this
analysis is developed for FE discretizations on
triangular grids and in~\cite{Vankasubmitted}, a general LFA technique
on simplicial grids for overlapping smoothers is presented.

\section{Numerical results}\label{sec:experiments}

In this section, we demonstrate the efficiency of the GMG method
for the MFD discretizations based on the
multiplicative Schwarz smoother and the inter-grid transfer operators
obtained from the modified N\'ed\'elec FE discretization. We also show
a local Fourier analysis for this kind of discretizations
to confirm the experimental results.

We first consider problem~\eqref{rot_rot_problem} with
$\kappa = 1$ on an equilateral triangular domain of unit
side-length. In Table~\ref{correspondence_equilateral}, we display the
smoothing factor $\mu$ and the two-grid $\rho_{2g}$ convergence factors
predicted by the LFA, together with the asymptotic
convergence factor computed by using a $W$-cycle on a target fine-grid
obtained after $10$ refinement levels. Since in practice it is worth
to know if we can use $V$-cycles instead of $W$-cycles, due to the high
computational cost of the latter, we also show the three-grid
convergence factors $\rho_{3g}$ predicted by LFA for $V$-cycles,
together with the asymptotic convergence factors experimentally
obtained. These results are presented for different numbers of
smoothing steps $\nu$.
\begin{table}
\begin{center}
\begin{tabular}{|c|c|c|c|c|c|} \hline & &
\multicolumn{2}{|c|}{$W$-cycle} & \multicolumn{2}{|c|}{$V$-cycle}\\ \hline
$\nu$ & $\mu^{\nu}$ & $\rho_{2g}$ & $\rho_h^W$ & $\rho_{3g}$ &
$\rho_h^V$ \\ \hline 1 & 0.462 & 0.331 & 0.330 & 0.337 & 0.334 \\
\hline 2 & 0.214 & 0.124 & 0.124 & 0.133 & 0.132 \\ \hline 3 & 0.099 &
0.070 & 0.069 & 0.072 & 0.071 \\ \hline 4 & 0.046 & 0.045 & 0.045 &
0.052 & 0.052 \\ \hline
\end{tabular}
\caption{Smoothing ($\mu$), two-grid ($\rho_{2g}$) and three-grid
($\rho_{3g}$) LFA convergence factors, together with measured $W$-cycle
and $V$-cycle asymptotic convergence rates, $\rho_h^W$ and $\rho_h^V$,
respectively, for an equilateral triangle and different numbers of
smoothing steps, $\nu$.}
\label{correspondence_equilateral}
\end{center}
\end{table} From the results in
Table~\ref{correspondence_equilateral}, we observe accurate
predictions of the asymptotic convergence factors. Moreover, an optimal behavior of $V$-cycle is shown,
since the obtained $V$-cycle convergence rates are very similar to the $W$-cycle convergence rates. We have seen that in this case, the
convergence factors are independent of how the smoothing
steps are distributed, and therefore we do not distinguish different
distributions in the table. Notice that very good convergence factors,
below $0.1$, are obtained by using a $V$-cycle with only three smoothing
steps.

\begin{figure}[htb] \centering
\includegraphics[scale = 0.45]{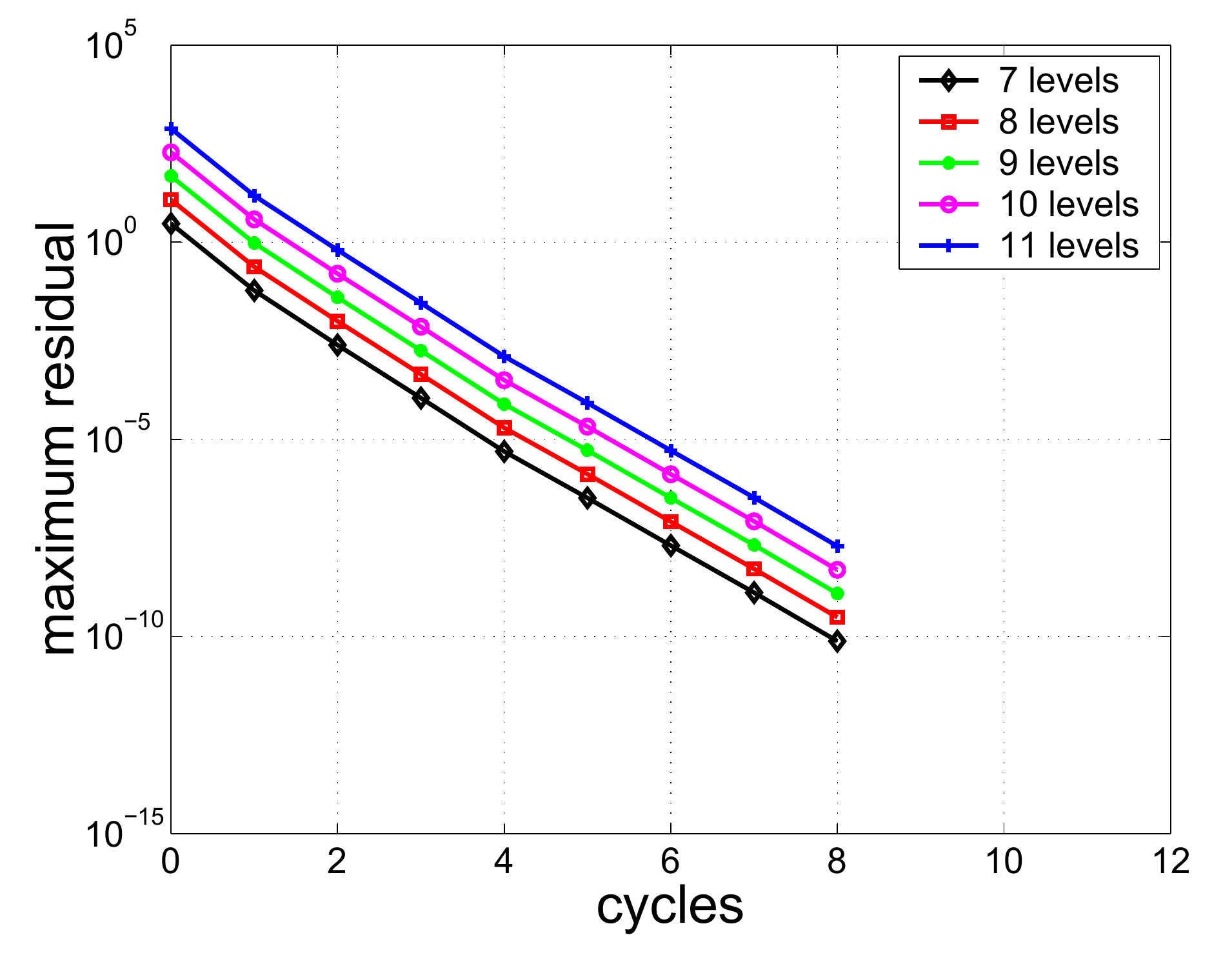}
\caption{History of the multigrid convergence of a $V(2,1)$-cycle for
different numbers of refinement levels.}
\label{h_convergence}
\end{figure}

Next, in Figure~\ref{h_convergence}, we display the history of the
MG convergence for different fine grids. We use a $V(2,1)$-cycle and the stopping criterion is
to reduce the initial residual by a factor of $10^{-10}$. As is a
well-known property of the MG methods, we observe the h-independent
convergence behavior.

\begin{table}[!htb]
\begin{center}
\begin{tabular}{|c|c|c|c|} \hline $\kappa$ & $\mu$ & $W$-cycle
($\rho_{3g}^W$) & $V$-cycle ($\rho_{3g}^V$) \\ \hline 1 & 0.099 & 0.070
& 0.072 \\ \hline $10^{-2}$ & 0.099 & 0.070 & 0.072 \\ \hline
$10^{-4}$ & 0.099 & 0.070 & 0.072 \\ \hline $10^{-6}$ & 0.099 & 0.070
& 0.072 \\ \hline $10^{-8}$ & 0.099 & 0.070 & 0.072 \\ \hline
\end{tabular}
\caption{Smoothing ($\mu$) and three-grid LFA convergence factors by
using $W$-cycle ($\rho_{3g}^{W}$) and $V$-cycle ($\rho_{3g}^{V}$) for an
equilateral triangular grid, with three smoothing steps ($\nu = 3$),
and different values of parameter $\kappa$.}
\label{robustness_k}
\end{center}
\end{table}
 To study the robustness of the proposed method with respect to
parameter $\kappa$, in Table~\ref{robustness_k}, we show the smoothing
and three-grid convergence factors, for both $W$- and $V$-cycle, predicted
by LFA for different values of $\kappa$ and by using three smoothing
steps. From the table, it is clear that the results are independent of $\kappa$.

\begin{figure}
\centering
 \includegraphics[width= 0.55\textwidth]{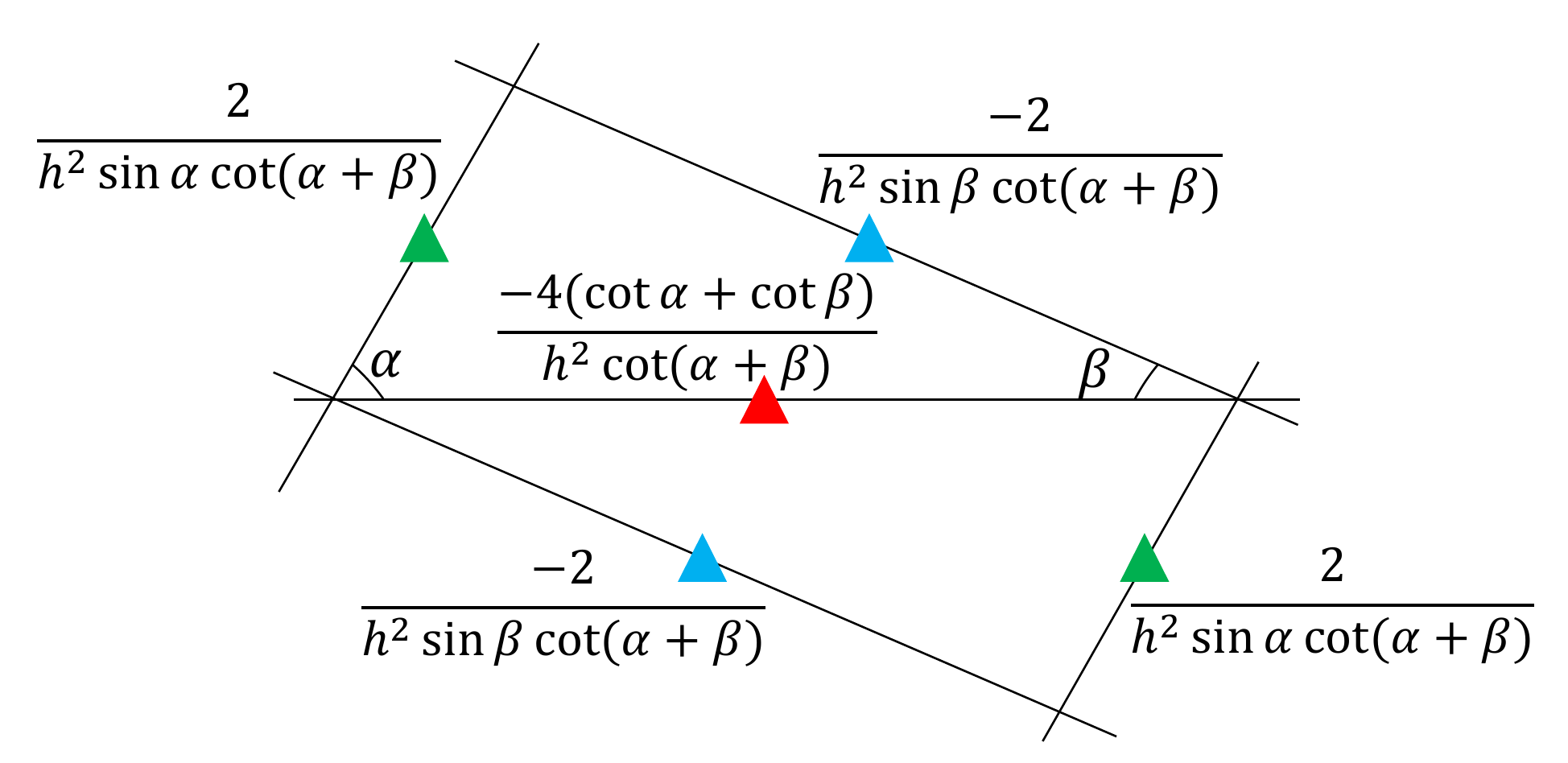}
 \caption{Stencil of \textbf{curl}-rot operator for an arbitrary
 triangulation characterized by angles $\alpha$ and $\beta$ by using
 mimetic finite differences.}
 \label{general_stencil}
 \end{figure}

 To show a more general applicability of the method, arbitrary
 structured triangular grids are considered.  These grids can be
 characterized by two angles $\alpha$ and $\beta$, and therefore, after
 simple computations we can obtain the stencil corresponding to
 \textbf{curl}-rot operator in terms of $\alpha$ and $\beta$, as we can see in
 Figure~\ref{general_stencil}.

\begin{figure}
\centering
 \includegraphics*[width=0.5\textwidth]{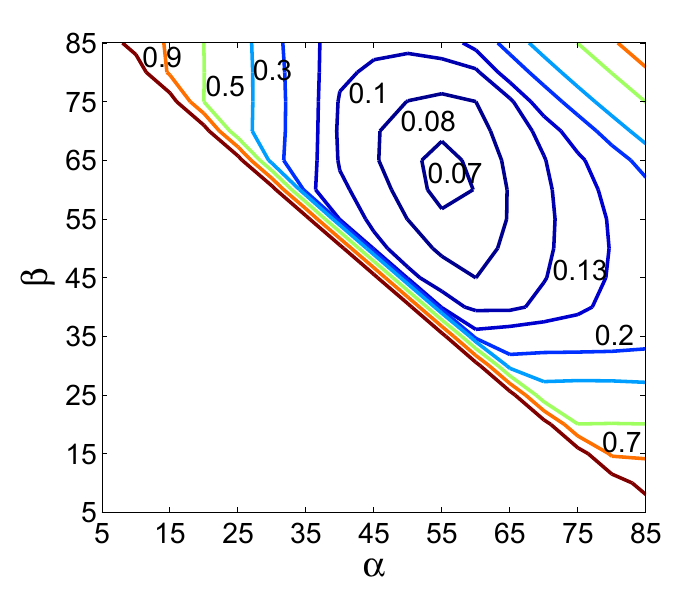}
 \caption{Three-grid convergence factors for $V$-cycle predicted by LFA
 for a wide range of triangles characterized by angles $\alpha$ and
 $\beta$.}
 \label{three_grid_angles}
 \end{figure}

 In this way, a systematic analysis with the LFA tool can
 be performed for a wide range of
 triangulations. In Figure~\ref{three_grid_angles}, we display
 the three-grid convergence factors predicted by the LFA for a wide range of triangular grids characterized by angles $\alpha$
 and $\beta$. For these results, a $V$-cycle with three smoothing steps
 has been considered.

 From Figure~\ref{three_grid_angles}, we observe a
 deterioration of the convergence factor when a small angle appears in
 the triangulation. This behavior is typical when point-wise smoothers
 are considered on a grid with anisotropy. It is possible
 to improve these convergence factors by using a relaxation parameter
 $\omega$. For example, if an isosceles triangular grid with angles
 $80^{\circ}$-$80^{\circ}$-$20^{\circ}$ is considered, we
 obtain a factor of $\rho_{3g}^V = 0.508$, but we can improve this
 result to $\rho_{3g}^V = 0.252$ by considering a relaxation parameter
 $\omega = 1.35$. These optimal parameters can be obtained for other
 triangulations by using the LFA. Of course, other well-known
 techniques can be applied to overcome the difficulties arising from
 the anisotropy of the grid, such as line-type smoothers, but this is
 not the focus of this work.

\section{Conclusions}\label{sec:conclusions}

In this work, we showed an equivalence between the MFD schemes on simplicial grids and some modified FE methods. This relation has been obtained for two model
problems in $\Hcurl$ and $\Hdiv$. This connection leads to immediate
convergence results for the MFD schemes
and also allows the construction of efficient multilevel
methods using the FE framework.  Based on the LFA, we theoretically predicted and numerically showed the robustness and the efficiency of the GMG method in
$\Hcurl$  for edge-based discretizations on simplicial grids.

\section*{Acknowledgements}
The work of Francisco J. Gaspar and Carmen Rodrigo is supported in
part by the Spanish project FEDER /MCYT MTM2013-40842-P and the DGA
(Grupo consolidado PDIE). The research of Ludmil Zikatanov is supported in
part by NSF DMS-1217142 and NSF DMS-1418843. Carmen Rodrigo gratefully
acknowledges the hospitality of the Department of Mathematics of The
Pennsylvania State University, where this research was partly carried
out.

\bibliographystyle{elsarticle-num-names}
\section*{\refname}
\bibliography{bib_rot_rot}



\end{document}